\documentclass[12pt]{amsart}
\usepackage{amsfonts}
\usepackage{amssymb}
\theoremstyle{definition}

\theoremstyle{remark}

\newcommand{\R}{\mathbb R}
\newcommand{\C}{\mathbb C}

\newcommand{\X}{\mathfrak X}

\newcommand{\ds}{\displaystyle}

\begin{document}

\centerline{\large\bf ON TOTALLY REAL SUBMANIFOLDS}

\vspace{0.1in}

\centerline{\large\bf Ognian T. Kassabov}

\vspace{0.1in}
\centerline{\large\bf University of Sofia}

\vspace{0.1in}
\centerline{\large\bf Received March, 9 1985.}

\vspace{0.4in}
{\sl Many autors have investigated totally real submanifolds of K\"ahlerian manifolds.
In this paper we generalize some results in this direction obtained by
B.-Y. Chen, K. Ogiue, M. Kon and K. Yano in [1] [4] [5] [6]. In particular,
we study semiparallel totally real submanifolds. The notion of semiparallel 
submanifolds was introduced by J. Deprez in [2] [3] as extrinsic analogue
for semisymmetric Riemannian spaces.}

\vspace{0.4in}
\centerline{\bf 1. Preliminaries.} 

\vspace{0.1in}
Let $\widetilde M$ be a $2m$-dimensional K\"ahlerian manifold with Riemannian metric 
$g$, complex structure $J$ and covariant differentiation $\widetilde\nabla$. An
$n$-dimensional submani\-fold $M$ of $\widetilde M$ is said to be a totally real
submanifold of $\widetilde M$, if for each point $p\in M$ the inclusion
$JT_pM \subset T_pM^{\perp}$ holds. Then $n \le m$. For $X,\,Y \in \X(M)$ we
write the Gauss formula:

$$
	\widetilde\nabla_X Y=\nabla_X Y +\sigma(X,Y) \,,
$$

\vspace{0.1in}\noindent
where $\nabla$ is the covariant differentiation on $M$ and $\sigma$ is a 
normal-bundle-valued symmetric tensor field on $M$, called the second fundamental form 
of $M$. The mean curvature vector $H$ of $M$ is defined by $H=(1/n)\,{\rm tr}\sigma$.
If $H=0$, $M$ is said to be a minimal submanifold of $\widetilde M$. In particular, if
$\sigma=0$, $M$ is called a totally geodesic submanifold of $\widetilde M$. For
$\xi \in \X(M)^{\perp}$ the Weingarten formula is given by

$$
	\widetilde\nabla_X \xi=-A_{\xi}X +D_X \xi  \,,
$$ 

\vspace{0.1in}\noindent
where $-A_{\xi}X$ (resp. $D_X\xi$) denotes the tangential (resp. the normal) component
of $\widetilde\nabla_X \xi$. It is well known that $g(\sigma(X,Y),\xi)=g(A_\xi X,Y)$ and $D$
is the covariant differentiation in the normal bundle. If $D_X\xi=0$ for each
$X \in \X(M)$, the normal vector field $\xi$ is said to be parallel. For a normal
vector field $\xi$, let $J\xi=P\xi+f\xi$, where $P\xi$ (resp. $f\xi$) denotes the tangential
(resp. the normal) component of $J\xi$. Then $f$ is an endomorphism of the normal
bundle and $f^3+f=0$. So if $f$ does not vanishes, it defines an $f$-structure in 
the normal bundle. If $Df=0$, i.e. $D_X f\xi-fD_X\xi=0$ for all $X \in \X(M)$, 
$\xi \in \X(M)^{\perp}$, the $f$-structure in the normal bundle is said to be
parallel. In this case it is not difficult to find that

$$
	\sigma(X,Y)=JA_{JX}Y=JA_{JY}X \,,   \leqno (1)
$$
$$
	D_X JY=J\nabla_X Y \,,   \leqno (2)
$$
$$
	A_\xi = 0 \,,   \leqno (3)
$$

\vspace{0.1in}\noindent
for $X,\,Y \in \X(M)$, $\xi \perp \X(M)\oplus J\X(M)$, see [6, p. 46]. We note, that 
if $n=m$, $f$ vanishes and (1) and (2) hold good. Because of (1), the equation 
of Gauss can be written in the form

$$
	\{ \widetilde R(x,y)z \}^t = R(x,y)z - [A_{Jx},A_{Jy}]z \,,
$$

\vspace{0.1in}\noindent
for $x,\,y,\,z \in T_pM$, $p\in M$ where $\widetilde R$ (resp. $R$) is the curvature 
tensor for $\widetilde M$ (resp. $M$) and $\{ \widetilde R(x,y)z \}^t$ denotes the
tangential component of $ \widetilde R(x,y)z $. In particular, if $\widetilde M$ is
a complex space form $\widetilde M(\mu)$, i.e. a K\"ahlerian manifold  of constant
holomorphic sectional curvature $\mu$, we find

$$
	\frac{\mu}4 x \wedge y = R(x,y) - [A_{Jx},A_{Jy}]\,.  \leqno (4)
$$

\vspace{0.1in}\noindent
If $[A_\xi,A\eta]=0$ for all $\xi,\,\eta \in \X(M)^\perp$, $M$ is said to have
commutative second fundamental forms. Then from (3) and (4) we obtain

\vspace{0.2in}
{\bf Lemma 1.} [6, p. 57] Let $M$ be a totally real submanifold of a complex space
form $\widetilde M(\mu)$. If the $f$-structure in the normal bundle is parallel, then
$M$ is of constant curvature $\mu/4$ if and only if $M$ has commutative second
fundamental forms.

\vspace{0.2in}
For minimal submanifolds we have

\vspace{0.2in}
{\bf Lemma 2.} [6, p. 57] Let $M$ be a totally real minimal submanifold with commutative 
second fundamental forms of a K\"ahlerian manifold $\widetilde M$. If the $f$-structure 
in the normal bundle is parallel, then $M$ is totally geodesic.

\vspace{0.1in}
We note that because of (1)

$$
	M \ {\rm is\ minimal\ if\ and\ only\ if } \ 
	\sum_{i=1}^{n} A_{Je_i}e_i = 0     \leqno (5)
$$

\vspace{0.1in}\noindent
for any orthonormal basis $\{ e_1,\hdots,e_n \}$ of a tangent space of $M$.

Let $R^\perp$ be the curvature tensor of the normal connection, i.e.

$$
	R^\perp(X,Y)\xi = D_XD_Y\xi - D_YD_X\xi - D_{[X,Y]}\xi \,,
$$

\vspace{0.1in}\noindent
for $X,\,Y \in \X(M),\, \xi \in \X(M)^\perp$. Using (2) we find that

$$
	R^\perp(X,Y)JZ = JR(X,Y)Z \,.    \leqno (6)
$$

\vspace{0.1in}
Let $\overline\nabla$ denote the covariant differentiation with respect to the
connection of van der Waerden - Bortolotti. For example

$$
	(\overline\nabla_X\sigma)(Y,Z)=D_X\sigma(Y,Z)-\sigma(\nabla_XY,Z)-\sigma(Y,\nabla_XZ) \,.
$$

\vspace{0.1in}\noindent
If $\overline\nabla\sigma=0$, then $M$ is said to have parallel second fundamental form or 
to be a parallel submanifold. More generally $M$ is said to be a semiparallel submanifold
of $\widetilde M$, if $\overline R(X,Y).\sigma =0$ [2] [3], where

$$
	(\overline R(X,Y).\sigma)(Z,U)=
			(\overline\nabla_X (\overline\nabla_Y \sigma))(Z,U) 
			- (\overline\nabla_Y(\overline\nabla_X\sigma))(Z,U) \,, 
$$ 

\vspace{0.1in}\noindent
or equivalently

$$
	(\overline R(X,Y).\sigma)(Z,U)= R^\perp(X,Y)\sigma(Z,U)-
			\sigma(R(X,Y)Z,U)-\sigma(Z,R(X,Y)U) \,.
$$

\vspace{0.1in}\noindent
In Euclidean spaces submanifolds of this kind are considered by J. Deprez [2] [3].
Using (1) and (6) we derive

\vspace {0.2in}
{\bf Lemma 3.} Let $M$ be a totally real submanifold with parallel $f$-structure in the
normal bundle of a K\"ahlerian manifold $\widetilde M$. Then $M$ is semiparallel
if and only if

$$
	R(x,y)A_{Jz}u= A_{Jz}R(x,y)u+A_{Ju}R(x,y)z
$$

\vspace{0.1in}\noindent
for all $x,\,y,\,z,\,u \in T_pM$, $p\in M$.

\vspace{0.1in}
On the other hand, as a generalization of the submanifolds with parallel mean curvature 
vector, the submanifolds with semiparallel mean curvature vector are defined by 
$R^{\perp}(X,Y)H=0$. We note that the class of submanifolds with semipa\-rallel mean 
curvature vector includes also the semiparallel submanifolds.

\vspace{0.4in}
\centerline{\bf 2. Submanifolds of constant curvature.}

\vspace{0.2in}
{\bf Proposition 1.} Let $M$ be an $n$-dimensional ($n>1$) totally real submanifold
of constant curvature $c$, with parallel $f$-structure in the normal bundle of a
K\"ahlerian manifold $\widetilde M$. If the mean curvature vector $H$ of $M$ is
semiparallel, then $M$ is minimal or flat.

\vspace{0.1in}
{\bf Proof.} First we note that $H \in J\X(M)$, because of (1). Now using (6) we obtain:

$$
	0=R^{\perp}(X,Y)H=-JR(X,Y)JH \,,
$$

\vspace{0.1in}\noindent
i.e.:

$$
	R(X,Y)JH=0 \,.
$$

\vspace{0.1in}\noindent
Since $M$ is of constant curvature $c$, this implies that

$$
	c \{g(y,JH)x-g(x,JH)y \}=0   \leqno (7)
$$

\vspace{0.1in}\noindent
for all $x,\,y \in T_pM$. Let $c\ne 0$. Putting in (7) $y=(JH)_p$, $x\perp y$, $x\ne 0$
we obtain $H_p=0$. Hence $M$ is minimal.

\vspace{0.1in}
Now we prove the main result in this section. 

\vspace{0.2in}
{\bf Theorem 1.} Let $M$ be an $n$-dimensional $(n>1)$ totally real semiparallel
submanifold of constant curvature $c$ with parallel $f$-structure in the normal
bundle of a K\"ahlerian manifold $\widetilde M$. Then $M$ is flat, i.e. $c=0$
or $M$ is a totally geodesic submanifold of $\widetilde M$. 

\vspace{0.1in}
{\bf Proof.} Since $M$ is of constant curvature $c$, $R(x,y)=c\, x \wedge y$ holds
good. Hence Lemma 3 implies:

$$
	\begin{array}{r}\vspace{0.1in}
		c\{ g(y,A_{Jz}u)x-g(x,A_{Jz}u)y \} = c\{ g(y,z)A_{Ju}x -g(x,z)A_{Ju}y  \\
								+ g(y,u)A_{Jz}x -g(x,u)A_{Jz}y \}
	\end{array}  \leqno (8)
$$

\vspace{0.1in}\noindent
for all $x,\,y,\,z,\,u \in T_pM$, $p\in M$. Let $c\ne 0$ and $\{ e_1,\hdots,e_n \}$ be 
an orthonormal basis of $T_pM$. We put $x=u=e_i$ in (8) and we add for $i=1,\hdots,n$
using (5) and Proposition 1. The result is

$$
	(n+1)cA_{Jz}y=0 \,.
$$

\vspace{0.1in}\noindent
Hence the assertion follows, because of (3).

\vspace{0.1in}
Using Lemma 1 we obtain

\vspace{0.2in}
{\bf Corollary.} Let $M$ be an $n$-dimensional $(n>1)$ totally real semiparallel
submanifold of constant curvature $c$, with parallel $f$-structure in the normal
bundle of a complex space form $\widetilde M(\mu)$. Then $M$ is flat, i.e. $c=0$
or $M$ is a totally geodesic submanifold of $\widetilde M(\mu)$, i.e. $c=\mu/4$.

\vspace{0.1in}
For parallel minimal submanifolds this Corollary is proved in [1], see also
[6, p.61].

\vspace{0.4in}
\centerline{\bf 3. Minimal submanifolds and the sign of the scalar curvature.}

\vspace{0.1in}
Let $M$ be an $n$-dimensional totally real minimal semiparallel submanifold
with parallel $f$-structure in the normal bundle of a complex space form 
$\widetilde M(\mu)$. According to (4) and Lemma 3:

$$
	\begin{array}{r}\vspace{0.1cm}
		\ds \frac\mu4 \{g(y,A_{Jz}u)x-g(x,A_{Jz}u)y \} + [A_{Jx},A_{Jy}]A_{Jz}u  \\ \vspace{0.1cm}
		   \ds = \frac\mu4 \{ g(y,z)A_{Ju}x-g(x,z)A_{Ju}y  \\  \vspace{0.1cm}
		                     +g(y,u)A_{Jz}x-g(x,u)A_{Jz}y \}  \\  \vspace{0.1cm}
		                +A_{Ju}[A_{Jx},A_{Jy}]z+A_{Jz}[A_{Jx},A_{Jy}]u
	\end{array}
$$

\vspace{0.1in}\noindent
holds good and hence we derive the relation:

$$
	\begin{array}{r}\vspace{0.1cm}
		\ds \frac\mu4 \{g(y,A_{Jz}A_{Ju}v)g(x,w)-g(x,A_{Jz}A_{Ju}v)g(y,w) \}   \\ \vspace{0.1cm}
																				 + g([A_{Jx},A_{Jy}]A_{Jz}A_{Ju}v,w)  \\ \vspace{0.1cm}
		   \ds = \frac\mu4 \{ g(y,z)g(A_{Jx}A_{Ju}v,w)-g(x,z)g(A_{Jy}A_{Ju}v,w)\\  \vspace{0.1cm}
		                    +g(y,A_{Ju}v)g(A_{Jz}x,w)-g(x,A_{Ju}v)g(A_{Jz}y,w)\}  \\  \vspace{0.1cm}
		                    +g([A_{Jx},A_{Jy}]z,A_{Jw}A_{Ju}v)+g(A_{Jz}[A_{Jx},A_{Jy}]A_{Ju}v,w)  \\  
	\end{array}   \leqno(9)
$$

\vspace{0.1in}\noindent
for all $x,\,y,\,z,\,u,\,v,\,w \in T_pM$, $p\in M$. Let $\{ e_1,\hdots,e_n \}$ be 
an orthonormal basis of $T_pM$. In (9) we put $ x=w=e_i$, $y=u=e_j$, $z=v=e_k$ and we
add for $i,\,j,\,k=1,\hdots,n$; this gives:

$$
	\begin{array}{r}\vspace{0.1cm}
		\ds \frac\mu4 (n+1)\sum_{i=1}^n {\rm tr}\,A_{Je_i}^2 -2\sum_{i,j=1}^n {\rm tr}\,A_{Je_i}^2A_{Je_j}^2-
															\sum_{i,j=1}^n {\rm tr}\,A_{Je_i}A_{Je_j}^2A_{Je_i} \\ \vspace{0.1cm}
		   +\ds 2\sum_{i,j=1}^n {\rm tr}\,A_{Je_i}A_{Je_j}A_{Je_i}A_{Je_j} =0 \,.   
	\end{array}   \leqno(10)
$$

\vspace{0.1in}\noindent
On the other hand it is not difficult to find that 

$$
	\sum_{i,j=1}^n {\rm tr}\,A_{Je_i}^2A_{Je_j}^2= \sum_{i,j=1}^n {\rm tr}\,A_{Je_i}A_{Je_j}^2A_{Je_i} \,.
$$

\vspace{0.1in}\noindent
Applying this to (10) we obtain:

$$
	\frac\mu4 (n+1)\sum_{i=1}^n {\rm tr}\,A_{Je_i}^2 +
	\sum_{i,j=1}^n \{ {\rm tr}\,(A_{Je_i}A_{Je_j}-A_{Je_j}A_{Je_i})^2-{\rm tr}\,A_{Je_i}^2A_{Je_j}^2 \} =0 \,.
$$

\vspace{0.1in}\noindent
Now, just as Theorem 8.1 in [6, p. 69], we can prove the following

\vspace{0.2in}
{\bf Theorem 2.} Let $M$ be a totally real minimal semiparallel submanifold with
parallel $f$-structure in the normal bundle of a complex space form $\widetilde M(\mu)$.
If the square of the length of the second fundamental form is constant (or equivalently,
if $M$ has constant scalar curvature $\tau$), then $M$ is totally geodesic or $\tau \ge 0$.
Moreover, if $\tau=0$, then $M$ is flat.

\vspace{0.4in}

\centerline{\bf 4. Commutative second fundamental forms.}

\vspace{0.1in}
To begin this section we note, that if $M$ is a totally real submanifold with commutative
second fundamental forms and parallel $f$-structure in the normal bundle of a complex space
form $\widetilde M(\mu)$, then the equation of Gauss (4) reduces to

$$
	R(x,y)= \frac{\mu}4 x \wedge y  \,.   \leqno (11)
$$

\vspace{0.1in}
{\bf Proposition 2.} Let $M$ be an $n$-dimensional ($n>1$) totally real semiparallel
submanifold with commutative second fundamental forms of a complex space
form $\widetilde M(\mu)$. If the $f$-structure in the normal bundle is parallel, 
then $M$ is totally geodesic or flat.

\vspace{0.1in}
{\bf Proof.} According to (11), $M$ is of constant curvature $\mu/4$. Now the assertion
follows from Theorem 1.

\vspace{0.1in}
More generally, from (11), Lemma 2 and Proposition 1 we derive

\vspace{0.2in}
{\bf Proposition 3.} Let $M$ be an $n$-dimensional ($n>1$) totally real submanifold
with commutative second fundamental forms and parallel $f$-structure in the normal
bundle of a complex space form $\widetilde M(\mu)$. If the mean curvature vector $H$ 
of $M$ is semiparallel, then $M$ is totally geodesic or $\mu=0$.

\vspace{0.2in}
Propositions 2 and 3 generalize some results in [6, p. 62].

If $n=2$ and $M$ has parallel mean curvature vector, we can weaken the assumptions
of Proposition 3. Namely, we have

\vspace{0.2in}
{\bf Proposition 4.} Let $M$ be a totally real surface with commutative second 
fundamen\-tal forms and parallel $f$-structure in the normal bundle of a K\"ahler
manifold $\widetilde M$. If the mean curvature vector of $M$ is parallel, then
$M$ is flat or totally geodesic.

\vspace{0.1in}
This follows from Lemma 2 and the following

\vspace{0.2in}
{\bf Proposition 5.} Let $M$ be a totally real surface with parallel $f$-structure 
in the normal bundle of a K\"ahler manifold $\widetilde M$. If the mean curvature
vector $H$ of $M$ is parallel, then $M$ is flat or minimal.

\vspace{0.1in}
{\bf Proof.} Since $H$ is parallel, it has constant length. Let $H\ne 0$, i.e.
$M$ is not minimal. As in Proposition 1 we have $R(X,Y)JH=0$. Hence 
$R(X,JH,JH,X)=0$, and since $H$ does not vanish, then $M$ is flat.

\vspace{0.2in}
{\bf Theorem 3.} Let $M$ be an $n$-dimensional ($n>1$) complete totally real submanifold
with parallel mean curvature vector and commutative second funda\-mental forms
of a $2m$-dimensional simply connected complete complex space form $\widetilde M(\mu)$.
If the $f$-structure in the normal bundle is parallel and $M$ is not totally geodesic,
then $M$ is a pythagorean product of the form

$$
	S^1(r_1)\times \hdots \times S^1(r_p)\times \R^{n-p} 
$$

\vspace{0.1in}\noindent
in a $\C^n$ in $\C^m$, where $1\le p\le n$.

\vspace{0.1in}
{\bf Proof.} According to Proposition 3, $\mu=0$. Then  $\widetilde M(\mu)$ is
(isometric to) $\C^m$ and the assertion follows from Theorem 7.1 in [6, p. 65].

\vspace{0.1in}
For parallel submanifolds, Theorem 3 is proved in [5], (see also [6, p. 66]).

\vspace{0.2in}
{\bf Corollary.} Under the same assumptions as in Theorem 3, if $M$ is compact,
it is pythagorean product of the form

$$
	S^1(r_1)\times \hdots \times S^1(r_n)
$$

\vspace{0.1in}\noindent
in a $\C^n$ in $\C^m$.

\vspace{0.4in}
\centerline{\bf References}

\vspace{0.1in}
\noindent
[1]. B.-Y, Chen and K. Ogiue. {\it On totally real submanifolds}. Trans. Amer.
Math. Soc., 193\,(1974), 257-266.

\noindent
[2]. J. Deprez. {\it Semi-parallel hypersurfaces}. Rend. Sem. Mat. Univer. Politec. Torino,
44\,(1986), 303-316. 

\noindent
[3]. J. Deprez. {\it Semi-parallel surfaces in Euclidean space} J. Geom. 25\,(1985), 192-200.

\noindent
[4]. M. Kon. {\it Totally real minimal submanifolds with parallel second fundamental
forms.} Atti Accad. Naz. Lincei, 57\,(1974), 187-189.

\noindent
[5]. K. Yano and M. Kon. {\it Totally real submanifolds of complex space forms II.}
K\"odai Math. Sem. Rep., 27\,(1976), 385-399.

\noindent
[6]. K. Yano and M. Kon. {\it Anti-invariant submanifolds.} Marcel Dekker, INC.
New York and Basel, 1976.

\vspace {0.7cm}
University of Sofia

Faculty of Mathematics and Mechanics

5, Anton Ivanov Street

1126 Sofia

BULGARIA

\end{document}